\documentclass{svproc}
\usepackage{amsfonts, amssymb, latexsym, amsmath, lscape, amscd, color, epsfig, graphicx, indentfirst, enumerate, amsxtra, amstext,todonotes}
\usepackage{todonotes}

\usepackage{url}

\begin{document}
\mainmatter              
\title{Operator splitting for semi-explicit differential-algebraic equations and port-Hamiltonian DAEs}
\titlerunning{Operator Splitting for PHS}  
%
\author{Andreas Bartel \and Malak Diab \and Andreas Frommer \and Michael Günther}
\authorrunning{A. Bartel et al.} 
%
 \tocauthor{Andreas Bartel, Malak Diab, Andreas Frommer, Michael Günther}
\institute{Institute of Mathematical Modeling, Analysis and Computational Mathematics, Bergische Universität Wuppertal, Gaußstraße 20, D-42119 Wuppertal, Germany\\
\email{[bartel,mdiab,frommer,guenther]@uni-wuppertal.de}}

\maketitle              

\begin{abstract}
Operator splitting methods allow to split the operator describing a complex dynamical system into a sequence of simpler sub-systems and treat each part independently. In the modeling of dynamical problems, systems of (possibly coupled) differential-algebraic equations (DAEs) arise. This motivates the application of operator splittings which are aware of the various structural forms of DAEs. Here, we present an approach for the splitting of coupled index-1 DAE as well as for the splitting of port-Hamiltonian DAEs, taking advantage of the energy-conservative and energy-dissipative parts. We provide numerical examples illustrating our second-order convergence results.
\keywords{operator splitting, differential algebraic equations, port-\-Ha\-milto\-nian systems}
\end{abstract}
\section{Introduction}
Operator splitting \cite{OSM_ODE} is one of the powerful numerical tools in solving dynamical systems. By splitting the problem into a sequence of subproblems, one can exploit the different structural or physical properties of each subsystem independently. The operator splitting method is often used for solving initial value problems in ordinary differential equations (ODEs). It combines numerical integration for subsystems into an efficient scheme for the overall problem. More precisely, for an ODE-IVP $$x'(t) = f (x,t) = f_1(x,t) + f_2(x,t), \quad x(t_0)=x_0,$$ the so-called Lie-Trotter splitting method sequentially solves dynamical systems driven by $f_1$ and $f_2$ with step sizes $h$, where the subsystems are connected via initial condition. Strang splitting \cite{OSM_ODE} solves the ODE system driven by $f_1$, $f_2$ and again $f_1$ with time steps $h/2$, $h$ and $h/2$ respectively.

The question of extending the method from ODEs to DAEs has been addressed for DAE systems in a decoupled form \cite{OSM_DAE}. However, in a more general setting, the DAE model is not given in a decoupled structure. Decoupling the DAE system is possible in principle but is not preferable for the application of the splitting schemes, since we want to preserve the structure of the system and exploit it for the subsystems. This paper examines splitting methods for general coupled index-$1$ DAEs. In addition, it considers DAE systems in the port-Hamiltonian framework.

\section{Coupled Semi-Explicit Index-$1$ DAEs}
Semi-explicit index-$1$ DAEs arise in particular in network modeling, where the coupling of two or more networks can occur
via algebraic and differential variables. Here we focus on the coupling of two networks modeled by an index-$1$ coupled system of DAEs of the form 
\begin{subequations}
\label{overall_system}
\begin{align}
    y'_1 &= f_1(y_1, y_2, z_1,z_2),  \label{eq:sub1a}\\ 
    0 &= g_1(y_1,y_2,z_1,z_2), \label{eq:sub1b}\\ 
    y'_2 &= f_2(y_1, y_2,z_1, z_2), \label{eq:sub2a}\\
    0 &= g_2(y_1,y_2,z_1,z_2), \label{eq:sub2b}
\end{align}
\end{subequations}
with
$$
\frac{\partial (g_1,g_2)}{\partial (z_1,z_2)}
$$ 
regular in a neighborhood of the solution. System \eqref{overall_system} equipped with initial values 
\begin{equation}
\label{cdt_initial}
    y(0) = y_0 =\begin{bmatrix} y_{1,0} \\y_{2,0} \end{bmatrix}
    \in \mathbb{R}^{n_y}, 
    \quad z(0) = z_0 = \begin{bmatrix} 
    z_{1,0} \\ z_{2,0} \end{bmatrix} \in \mathbb{R}^{n_z}.
\end{equation}
is assumed to have a unique solution $y=(y_1,y_2)^\top : [0,T] \rightarrow \mathbb{R}^{n_y}$, $z=(z_1,z_2)^\top : [0,T] \rightarrow \mathbb{R}^{n_z}$ where $[0,T]=:\mathcal{I} $ is a finite interval. The functions $f_1$, $f_2$, $g_1$ and $g_2$ are supposed to be sufficiently differentiable in the neighbourhood of the solution. To apply a splitting technique, we introduce the following ``doubled''  (the algebraic constraints appear twice) decomposition of the DAE overall system \eqref{overall_system} 
\begin{equation}
    \label{DAE_decomposition}
    \begin{pmatrix}
        y'_1 \\ 0 \\ y'_2 \\ 0 
    \end{pmatrix} = \begin{pmatrix}
        f_1 \\ 2g_1 \\ f_2 \\ 2g_2 
    \end{pmatrix} = \underbrace{\begin{pmatrix}
        f_1 \\ g_1 \\ 0 \\ g_2
    \end{pmatrix}}_{\text{subsystem 1}} + \underbrace{\begin{pmatrix}
        0 \\ g_1 \\ f_2 \\ g_2
    \end{pmatrix}}_{\text{subsystem 2}}.
\end{equation}
Based on the index-$1$ condition for the overall system and using the implicit function theorem, there exists a continuously differentiable function $\varphi: \mathbb{R}^{n_y} \rightarrow \mathbb{R}^{n_z}$ such that $(z_1,z_2) = \varphi(y_1,y_2)$ which consequently allows to reduce the subsystems 
in \eqref{DAE_decomposition} to the ODEs
\begin{equation}
  \label{equivODE}
    \begin{cases}
        y'_1 = f_1(y_1,y_2,\varphi_1(y_1,y_2)) \\
        y'_2 = 0
    \end{cases}
    \quad \text{and} \qquad
    \begin{cases}
        y'_1 = 0 \\
        y'_2 = f_2(y_1,y_2,\varphi_2(y_1,y_2))
    \end{cases}.
\end{equation}
The structure of the systems \eqref{equivODE} leads to an operator splitting method based on the doubled decomposition \eqref{DAE_decomposition}, where the subsystems are given by
\begin{align}
       \label{subsystems}
       y'_1 &= f_1(y_1, y_2, z_1,z_2), & \quad y'_1 &= 0, \nonumber\\ 
       0 &= g_1(y_1,y_2,z_1,z_2), & \quad 0 &= g_1(y_1,y_2,z_1,z_2), \\ 
       y'_2 &= 0, & \quad y'_2 &= f_2(y_1, y_2,z_1,z_2), \nonumber\\
       0 &= g_2(y_1,y_2,z_1,z_2), & \quad 0 &= g_2(y_1,y_2,z_1,z_2), \nonumber
\end{align}
which is in this case equivalent to an ODE operator splitting \cite{OSM_ODE} and therefore similar convergence outcomes will result. We note that with this splitting the algebraic constraints are solved twice. For efficiency reasons, this is not a major limitation if the dimension of the algebraic constraints is small compared to the dimension of the dynamic quantities. 


\begin{proposition}
\label{prop:convergence-dedicated-coupling} 
Solving the two subsystems of \eqref{subsystems} sequentially on subintervals $[t_n,t_{n+1}]$ of $\mathcal{I}$ with at least first order convergent time integration methods yields a first order convergent method on $\mathcal{I}$. Furthermore, applying Strang splitting with at least second order time integration schemes yields a method of second order. 
\end{proposition}

Proposition \ref{prop:convergence-dedicated-coupling} follows directly from known convergence results for ODE operator splitting; see \cite{OSM_ODE} and \cite{OSM_DAE}.

\begin{remark}
    If we have only coupling via differential variables, i.e., 
    $\frac{\partial f_1}{\partial z_2}= \frac{\partial g_1}{\partial z_2}=0,\,\, \text{and} \,\,
    \frac{\partial f_2}{\partial z_1}= \frac{\partial g_2}{\partial z_1}=0$ in system \eqref{overall_system}, then we can proceed with the splitting without considering the doubled decomposition. In 
    this case, the subsystems in \eqref{subsystems} are given as follows: 
    \begin{align}
       y'_1 &= f_1(y_1, y_2, z_1), & \quad y'_1 &= 0, \nonumber\\ 
       0 &= g_1(y_1,y_2,z_1), & \quad y'_2 &= f_2(y_1, y_2,z_2), \\ 
       y'_2 &= 0, & \quad 0 &= g_2(y_1,y_2,z_2). \nonumber
\end{align}
\end{remark}
\begin{remark}
    A particular coupling type of system \eqref{overall_system}, where we have a dedicated coupling equation, has the form
    \begin{subequations}
      \label{overall_system2}
       \begin{align}
            y'_1 &= f_1(y_1, y_2, z_1),  \\ 
            0 &= g_1(y_1,y_2,z_1,u), \\ 
            y'_2 &= f_2(y_1, y_2, z_2), \\
            0 &= g_2(y_1,y_2,z_2,u), \\ 
            0 &= k(y_1,y_2,z_1,z_2), 
        \end{align}
    \end{subequations}
    where $k$ is a coupling equation and $u$ a dedicated Lagrangian coupling variable. Such coupled systems arise, for example, in circuit simulation or multibody dynamics~\cite{argu01}, and an operator splitting scheme can be similarly applied by considering the decomposition
    \begin{equation}
    \begin{pmatrix}
        y'_1 \\ 0 \\ y'_2 \\ 0 \\ 0
    \end{pmatrix} = \begin{pmatrix}
        f_1 \\ 2g_1 \\ f_2 \\ 2g_2 \\ 2k
    \end{pmatrix} = \underbrace{\begin{pmatrix}
        f_1 \\ g_1 \\ 0 \\ g_2 \\ k
    \end{pmatrix}}_{\text{subsystem 1}} + \underbrace{\begin{pmatrix}
        0 \\ g_1 \\ f_2 \\ g_2 \\ k
    \end{pmatrix}}_{\text{subsystem 2}} .
\end{equation}
\end{remark}

\section{Port-Hamiltonian DAEs}
In this section we are interested in applying the operator splitting method to  port-Hamiltonian systems (PHS) of differential-algebraic equations. An operator splitting for PHS-ODEs is considered in \cite{OSM2023} based on the natural decomposition given by the power-conserving part and the dissipative part. Unlike ODEs, DAEs suffer from singularity due to the presence of possibly hidden constraints. Different approaches like index reduction and decoupling the system---after which the system structure might be lost---allow for the application of ODE numerical methods to DAEs. Our goal is to avoid decoupling, so that we can preserve the system's structure.
\begin{definition}
    For $t \in \mathcal{I}$ and a state variable $x: \mathcal{I} \rightarrow \mathbb{R}^n$, a PHS with  quadratic Hamiltonian $H(x)=\frac{1}{2} x^\top Q^\top E x$ is given by
    \begin{subequations}
      \label{PHS}
      \begin{align}
        Ex'(t) &= (J-R)Qx(t) + Bu(t), \\
        y(t) &= B^\top Qx(t)
      \end{align}
    \end{subequations}
    where $E,Q,J$ and $R$ are real $n \times n$ matrices satisfying $E^\top Q = Q^\top E$, $R=R^\top \geq 0$, $J=-J^\top$, $B$ is $n \times m$ matrix and $Q$ is assumed to be symmetric positive definite. The functions 
    $u$ and $y$ are referred to as the inputs and outputs,  respectively. $E$ might not have a full rank.
\end{definition}

A fundamental property for the solution of \eqref{PHS} is the dissipativity inequality
\begin{equation}
    \label{dissipativity_ineq}
    \frac{d}{dt} H(x(t)) = \partial_x H(x(t))x'(t) = -x^\top Q^\top RQx + y^\top (t)u(t) \leq y^\top (t)u(t).
\end{equation}
The dissipativity inequality \eqref{dissipativity_ineq} can be also written in integral form
\begin{equation}
    \label{dissipativity_integral}
    H (x(t_0+h)) - H (x(t_0)) \leq \int_{t_0}^{t_0+h} y(s)^\top u(s) \, ds. 
\end{equation}

For the splitting of a system of the form \eqref{PHS}, we note that the regularity of the matrix pencils $\{E,J\}$ and $\{E,R\}$ is, as encountered in applications, not guaranteed. As an alternative, we adapt the idea of $\varepsilon$-embedding methods \cite{IVP_DAE95} and apply the splitting method to the perturbed system
\begin{subequations}
      \label{PHS_reg}
      \begin{align}
        E_\varepsilon x'(t) &= (J-R)Qx(t) + Bu(t),  \\
        y(t) &= B^\top Qx(t)
      \end{align}
\end{subequations}
with the perturbed flow matrix $E_\varepsilon = E+\varepsilon I$ representing a regularization of $E$. The system of equations \eqref{PHS_reg} is a PHS-ODE with the Hamiltonian $H_\varepsilon (x)=\frac{1}{2}x^\top Q^\top E_\varepsilon x$. Denoting by the pair $(x_\varepsilon , y_\varepsilon)$ the solution of system \eqref{PHS_reg}, we have
\begin{eqnarray*}
    H(
         x_\varepsilon(t_0&+&h)) 
         - H(x_\varepsilon(t_0)) = 
         \Big(   H(x_\varepsilon(t_0+h)) -  H_\varepsilon (x_\varepsilon(t_0+h)) \Big) 
         \\ 
         & &  \mbox{} + \Big( H_\varepsilon (x_\varepsilon(t_0+h))  - H_\varepsilon (x_\varepsilon(t_0)) \Big)  
         + \Big( H_\varepsilon (x_\varepsilon(t_0)) - H(x_\varepsilon(t_0)) \Big) 
         \\ 
         & \leq &
         - \tfrac{\varepsilon}{2} x_\varepsilon(t_0+h)^\top Q x_\varepsilon(t_0+h) + \int_{t_0}^{t_0+h} \!\!\! y_\varepsilon(s)^\top u(s) \, ds 
         + \tfrac{\varepsilon}{2} x_\varepsilon(t_0)^\top Q x_\varepsilon(t_0) 
         \\ 
         &\leq &
         \int_{t_0}^{t_0+h} y_\varepsilon(s)^\top u(s) \, ds + 
         \frac{\varepsilon}{2} \Big|  \| x_\varepsilon(t_0+h) \|_Q -  \| x_\varepsilon(t_0) \|_Q\Big|, 
\end{eqnarray*}
where $\| \cdot \|_Q$ is the $Q$-norm. We point out that the 
second term on the right-hand side, which is of order $\epsilon \cdot h$,
is negligible as $\varepsilon \rightarrow 0$. The solution $x_\varepsilon$ of the perturbed system \eqref{PHS_reg} converges to that of system \eqref{PHS}, see \cite{IVP_DAE95}, and thus $y_{\varepsilon}$ converges to $y$ for $\varepsilon \rightarrow 0$. Furthermore, the energy function $y_{\varepsilon}^\top \cdot u$ cannot increase indefinitely, therefore using the dominated convergence theorem, $\int_0^{h} y_{\varepsilon} (s)^\top u(s) \, ds$ converges to $\int_0^{h} y (s)^\top u(s) \, ds$ for $\varepsilon  \rightarrow 0$. So we may say that the dissipation inequality  
\eqref{dissipativity_ineq} is satisfied in the limit and that it is increasingly less violated as $\varepsilon \to 0$.

The splitting of the port-Hamiltonian system \eqref{PHS_reg} is performed by splitting the right-hand side into an energy-preserving part $f_1(x,t):=JQx(t)$ and a dissipative part $f_2(x,t):=-RQx(t)+Bu(t)$. The Strang splitting method yields the approximate value $w_{\varepsilon}(h/2)$ and consequently offers a second-order approximation for the exact solution, summarized as follows:
\begin{align*}
        z'_\varepsilon (t) &= f_1(z_\varepsilon,t), \,\, y_{z,\varepsilon} = B^\top Qz(t), \,\,\,\, z_0 = x_0, \,\,\, t\in[0,\tfrac{h}{2}], \\ 
        v'_\varepsilon (t) &= f_2(v,t), \,\, y_{v,\varepsilon}
        = B^\top Qv(t), \,\,\,\, v_0 = z(\tfrac{h}{2}), \,\, t \in [0,h], \\ 
        w'_\varepsilon(t) &= f_1(w,t), \,\, y_{w,\varepsilon} = B^\top Qw(t), \,\, w_0 = v(h), \,\,\, t\in[0,\tfrac{h}{2}]. 
\end{align*}  

If the port-Hamiltonian system \eqref{PHS_reg} is solved numerically using, for instance, Strang splitting, provided that the used numerical integration method is at least of second order, then the dissipativity inequality is preserved on a discrete level \cite{OSM2023}. For the exact solution it is given by
\begin{eqnarray*}
   H_\varepsilon(w_\varepsilon(\tfrac{h}{2})) - H_\varepsilon(x_0) 
    &= & \big( H_\varepsilon(w_\varepsilon(\tfrac{h}{2})) - H_\varepsilon (w_\varepsilon(0)) \big) + \big( H_\varepsilon(v_\varepsilon(h)) - H_\varepsilon(v_\varepsilon(0)) \big)   
    \\ 
    & &  \mbox{} 
   + \big( H_\varepsilon(z_\varepsilon(\tfrac{h}{2})) - H_\varepsilon(x_0) \big) 
   \\ 
  &\leq& 
     \int_0^{h} y_{v,\varepsilon}(s)^\top u(s) \, ds  .
\end{eqnarray*}
Note that only the second step in the Strang splitting method contributes to the inequality, the first and third being energy-conserving.

\section{Numerical Results}
We show numerical results demonstrating the order of convergence for the operator splitting method considered in the previous sections.
\begin{example}[{{\rm Coupled LC Oscillator}}]
    We consider the coupled, linear circuit in Figure \ref{fig:circuit_example1} consisting of two resistors $R_1$, $R_2 > 0$, two capacitors $C_1$, $C_2 >0$ and two inductors $L_1$, $L_2 >0$. By $u_1,u_2,u_3$ and $u_4$ we denote the node potentials, $j_{L,1}, j_{L,2}$ are the currents through $L_1$ and $L_2$ respectively, and $j_{co}$ is a coupling current. 
    The circuit is modeled by the following two subsystems.
    \begin{align*}
      \text{subsystem1: } & &  \text{subsystem 2:} \\
       C_1\frac{d}{dt}u_1 - \frac{1}{R_1} (u_2 - u_1) &= 0, & \quad  -\frac{1}{R_2} (u_4 - u_3) + j_{L,2} - j_{co} &= 0, \\
       \frac{1}{R_1} (u_2 - u_1) + j_{L,1} + j_{co} &= 0, & \quad  L_2 \frac{d}{dt} j_{L,2} - u_3 &= 0, \\
       L_1 \frac{d}{dt} j_{L,1} - u_2 &= 0,  & \quad  C_2\frac{d}{dt}u_4 + \frac{1}{R_2} (u_4 - u_3) &= 0,  \\
        & & \quad  u_2 - u_3 &= 0.
    \end{align*}
    Both subsystems, as well as the overall coupled system are of index one. We compute a reference solution for the overall system using the midpoint rule with time step $h=10^{-7}$. A comparison is then done by solving sequentially subsystems $1$ and $2$ using the Lie-Trotter and Strang splitting approaches. The time integration of the subsystems is performed using the midpoint rule as a second order scheme, with different choices for the time steps $h$ between $10^{-5}$ and $10^{-7}$. The results given in Figure~\ref{fig:order_sol_example1} verify that the orders of convergence of the Lie-Trotter and Strang splitting approaches are indeed $1$ and $2$, respectively.   
    \begin{figure}[ht]
        \centering
        \includegraphics[scale=0.28]{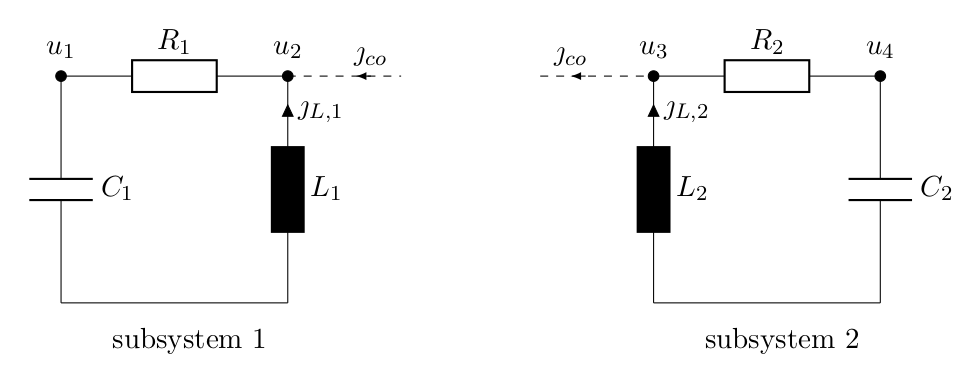}
        \caption{Subsystems of coupled oscillator circuit via coupling equation.}
        \label{fig:circuit_example1}
    \end{figure} 
    \begin{figure}[ht]
        \centering
        \includegraphics[scale=0.4]{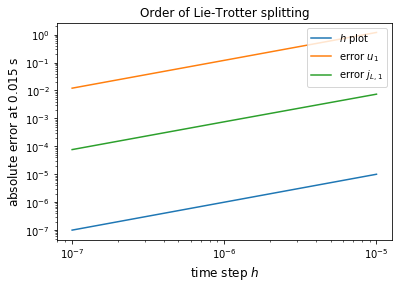}
        \includegraphics[scale=0.4]{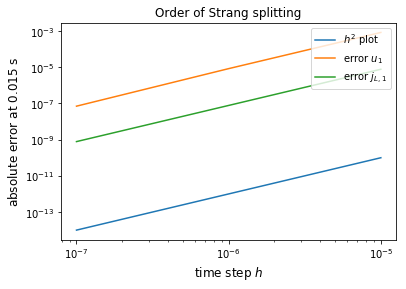}
        \caption{Circuit \ref{fig:circuit_example1}: order of convergence for some of the state variables.}
        \label{fig:order_sol_example1}
    \end{figure}
\end{example}
\begin{example} \label{ex:circuit2}
    For $x^\top = (u_1,u_2,u_3,j_1,u_4,u_5)$ the electric circuit in Figure \ref{fig:circuit_example3} is modeled  by 
    \[
    Ex' = (J-R)x+Bi(t) \enspace \text{with } E=\text{diag}(0,C_1,0,L_1,0,C_2).
    \]
    \begin{figure}
        \centering
        \includegraphics[scale=0.28]{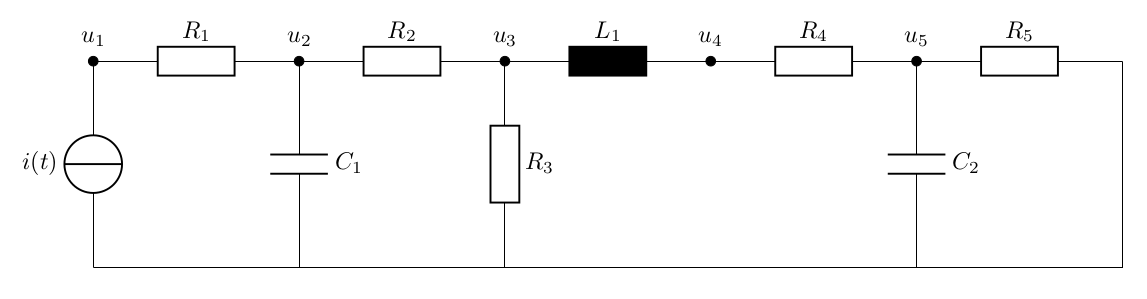}
        \caption{Electric circuit for Example~\ref{ex:circuit2}}
        \label{fig:circuit_example3}
    \end{figure}
\end{example}

\begin{figure}[ht]
        \centering
        \includegraphics[scale=0.4]{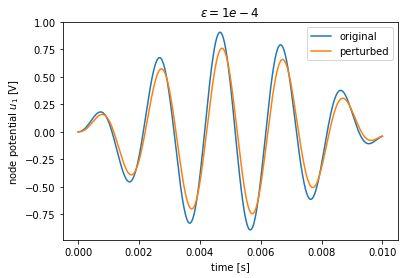}
        \includegraphics[scale=0.4]{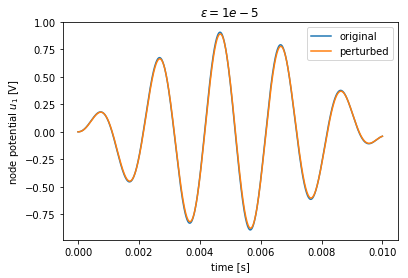}
        \caption{Plots for $u_1$ obtained by solving the original and perturbed systems.}
        \label{fig:eps=1e-4}
\end{figure}

\begin{figure}[ht]
    \centering
    \includegraphics[scale=0.4]{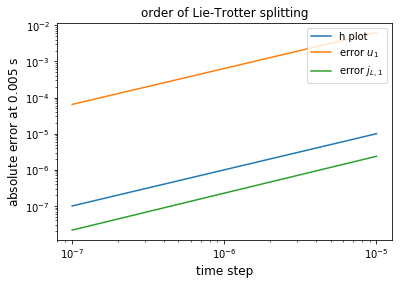}
    \includegraphics[scale=0.4]{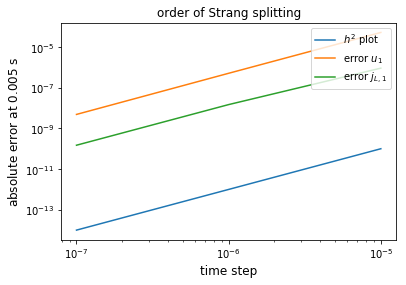}
    \caption{Order of convergence for some of the state variables of the circuit in Figure \ref{fig:circuit_example3}}
    \label{fig:circuit3_conv}
\end{figure}

Using $R_1=R_2=R_3=R_4=0.5, R_5=5, C_1=C_2=5\cdot 10^{-4}, L=20$ and
$i(t)=\sin(2\pi \cdot 50t) \cdot \sin(2\pi \cdot 500t)$, we compute a reference solution for the model equations describing the circuit in Fig.~\ref{fig:circuit_example3} using the mid-point rule and a step size of $h=10^{-8}$. In the first step, we compare the solution of the original system to the solution of the perturbed PHS-DAE problem \eqref{PHS_reg} for small values of $\varepsilon$. We solve using the implicit Euler method for $h=10^{-5}$, plots are shown in Fig.~\ref{fig:eps=1e-4} for $\varepsilon = 10^{-4}$ and $\varepsilon = 10^{-5}$ respectively. We observe the convergence of the solution of the perturbed problem to the solution of the original problem as $\varepsilon$ approaches zero. Similar to the previous example, we implement the Lie-Trotter and Strang splittings. For time integration we use the midpoint rule with $h=10^{-5}$, and the convergence results are plotted in Fig.\ref{fig:circuit3_conv}, illustrating again the first and second order convergence of the two splitting methods.


\section*{Conclusions and Outlook} In this paper we introduced a ``doubled'' decomposition for semi-implicit DAEs, on which our operator splitting is based. By doing this, one is able to work with each smaller part more effectively, using techniques that are both efficient and that keep the essential characteristics of the original problem. We additionally applied the operator splitting to PHS-DAE after perturbation to avoid non-regularity. The numerical tests we considered demonstrate second-order convergence for the Strang splitting  approach.
Further work will consider the non-regularity of the split DAE matrix pencils, index-$2$ DAEs and higher order cases.

%
%

\end{document}